\documentclass{amsart}      

\let\labelindent\relax
\usepackage{graphicx}   
\usepackage{epsfig} 
\usepackage[T1]{fontenc} 
\usepackage{amssymb} 
\usepackage{caption}
\usepackage{subcaption}
\usepackage{booktabs}
\usepackage{amsmath} 
\usepackage[pdf]{pstricks}
\usepackage{algorithm}
\usepackage[noend]{algpseudocode}
\usepackage{graphicx}
\usepackage{lipsum}%
\graphicspath{{./Figures/}}
\usepackage{xcolor}
\usepackage{hyperref}
\usepackage{float}
\usepackage{enumitem}
\usepackage[backend=biber,style=ieee,sorting=none]{biblatex}
\newcommand{\RNum}[1]{\uppercase\expandafter{\romannumeral #1\relax}}

\usepackage{color}

\title{\LARGE \bf Long Horizon Risk-Averse Motion Planning: A Model-Predictive Approach}
\author{Chris van der Ploeg$^{1,2}$, Robin Smit$^{1}$, Arjan Teerhuis$^{1}$, Emilia Silvas$^{1,3}$} 
\thanks{$^{1}$Netherlands Organisation for Applied Scientific Research, Integrated Vehicle Safety Group, 5700 AT Helmond, The Netherlands.}
\thanks{$^{2}$Eindhoven University of Technology, Dynamics and Control Group, Mechanical Engineering Dept., P.O. Box 513, 5600 MB, Eindhoven, The Netherlands.(e-mail: \href{mailto:c.j.v.d.ploeg@tue.nl}{c.j.v.d.ploeg@tue.nl}).}
\thanks{$^{3}$Eindhoven University of Technology, Control Systems Technology Group, Mechanical Engineering Dept., P.O. Box 513, 5600 MB, Eindhoven, The Netherlands.}

\date{July 2022}
\begin{document}
\maketitle
\thispagestyle{empty}
\pagestyle{empty}
\begin{abstract}
    Developing safe automated vehicles that can be proactive, safe, and comfortable in mixed traffic requires improved planning methods that are risk-averse and that account for predictions of the motion of other road users.
    To consider these criteria, in this article, we propose a non-linear model-predictive trajectory generator scheme, which couples the longitudinal and lateral motion of the vehicle to steer the vehicle with minimal risk, while progressing towards the goal state. The proposed method takes into account the infrastructure, surrounding objects, and predictions of the objects' state through artificial potential-based risk fields included in the cost function of the model-predictive control (MPC) problem. This trajectory generator enables anticipatory maneuvers, i.e., mitigating risk far before any safety-critical intervention would be necessary. The method is proven in several case studies representing both highways- and urban situations. The results show the safe and efficient implementation of the MPC trajectory generator while proving its real-time applicability. 
\end{abstract}
\section{Introduction}
Increasing levels of vehicle automation are projected to reduce significantly the number of fatalities caused by collisions and maximize driver comfort and efficiency~\cite{Papadoulis201912}. To achieve this, it is of great importance for automated vehicles to negotiate and navigate in all traffic situations, dangerous or not, with minimal risk. While \textit{minimal risk} forms an essential argument in terms of safety,  \textit{progression} towards the desired destination forms a large additional factor for human acceptance of automated vehicles. For example, an excessively slow-driving vehicle, which is arguably less susceptible to risks, may severely reduce acceptance~\cite{DeFreitas2021}. While planning the motion of a vehicle, these factors of risk and progression need to be taken into account. 

Although various automated driving architectures exist, there is yet no expert consent on an optimal one. Nevertheless, in general, motion planning for an automated vehicle can be divided into global (e.g., \textit{route planning}) and local planning (e.g., \textit{trajectory generation})~\cite{González20161135}. The task of global planning, i.e., finding a route from some initial location to a destination, occurs on a large road segment level and is generally done through path search algorithms such as $\text{A}^*$~\cite{Likhachev2009933}. It is thereafter the task of a trajectory generator to process these global points and create driveable, safe, and comfortable trajectories to traverse the path towards the destination while satisfying the requirements of minimal risk and progression. 

Optimization-based methods, more specifically MPC algorithms, are often used to solve the trajectory generation problem. These methods explicitly take into account the vehicle model, as they can optimally balance multiple objectives and various constraints. Several works deploy an MPC-based trajectory generator, where the road boundaries and dynamical object avoidance are either formulated through hard constraints~\cite{Gutjahr20171586} or through the use of artificial potential fields (APFs)~\cite{Ji2017952}. Note, that the notion of \textit{risk} fits well with the control-theoretic notion of artificial potential fields, first introduced in~\cite{Khatib1985500}. These APFs allow the analytical formulation of repulsive fields, i.e., regions to stay away from, and attractive fields, i.e., regions to be attracted towards. Several works in literature advocate the use of these APFs as a replacement for hard anti-collision constraints, due to the potential non-convexity of the constraints and, hence, a higher computational burden~\cite{Dixit20181061}. Moreover, artificial potential fields circumvent the issue of the potential infeasibility of the optimization problem caused by unavoidable collisions, since they act in the problem cost function and, hence, a collision would not lead to the infeasibility of the problem. In this case, a minimal-risk form of damage mitigation would result in a solution, instead of no solution at all. 

Model-predictive APF-based trajectory generators have been proposed in several different contexts in literature, ranging from emergency collision avoidance~\cite{Ji2017952} up to longer-term planning, e.g., highway maneuvering~\cite{9366415,Huang2016232,Hang202014458}.  Many of these solutions, however, only consider relatively short horizons (e.g., 1.6s~\cite{Dixit20181061}) due to, e.g., the high fidelity of the used model. To make informed decisions and actions, a trajectory generator is expected to look ahead for a sufficient distance or time. For reference, a human driver would, depending on the context, look ahead and anticipate on the environment $3-4$ seconds up to tens of seconds~\cite{Calvert202042}. Another shortcoming in the literature is related to decision-making modules or pre-planned trajectories, which accompany the model-predictive trajectory generator, and are deployed to decide and give directions on, e.g., changing a lane~\cite{Hang202014458,Dixit20181061}. Involving separate components to pre-determine actions or directions of the trajectory generator may result in a loss of optimality (i.e., increased risk or reduced progression). 

\textit{Our contributions:} In summary, there exist various solutions for model-predictive minimal risk trajectory planning. However, there does not yet exist an approach that finds minimal risk decisions and actions in an optimal manner (i.e., without the use of additional finite state machines or pre-defined trajectories).
As such, the contributions of our work are focused on the following three points:
\begin{enumerate}[label=(\roman*)]
    \item {\textit{Risk-averse trajectory planning without discrete decisions or pre-defined trajectories:} we propose a monolithic trajectory planner which can undertake anticipatory complex risk-averse maneuvers towards a destination without the need for discrete decisions or pre-defined trajectories to follow. }
     \item {\textit{Multiple operational domains:} the planner is simulated both in relevant highway scenarios, as well as an urban VRU crossing scenario, incorporating both timely and late object detection (e.g., bad weather effects on sensors). This shows its applicability to a variety of operational domains and conditions.}
     \item {\textit{Real-time implementation:} the algorithm is implemented in the CasADi~\cite{Andersson2018} optimization framework using an interior point solver (IPOPT~\cite{Wachter200625}) and is deployed in real-time, showing the feasibility of solving such a problem in experimental settings.}
\end{enumerate} 
The outline of our work is as follows. First, in Section~\ref{sec:problemstatement}, a mathematical description of the problem statement is introduced. In Section~\ref{sec:method}, the optimization problem is formulated through the notions of risk and progression. In Section~\ref{sec:simulation}, a simulation study is presented in both highway and urban scenarios, showing the expected minimal-risk results of our scenarios and their robustness in the case of late object detections. Finally, Section~\ref{sec:conclusion} concludes the work.
\section{Preliminaries and Problem Statement}\label{sec:problemstatement}

\begin{figure}[h!]
    \centering
    \begin{subfigure}{\columnwidth}
    \centering
        \def\svgwidth{0.75\columnwidth}
        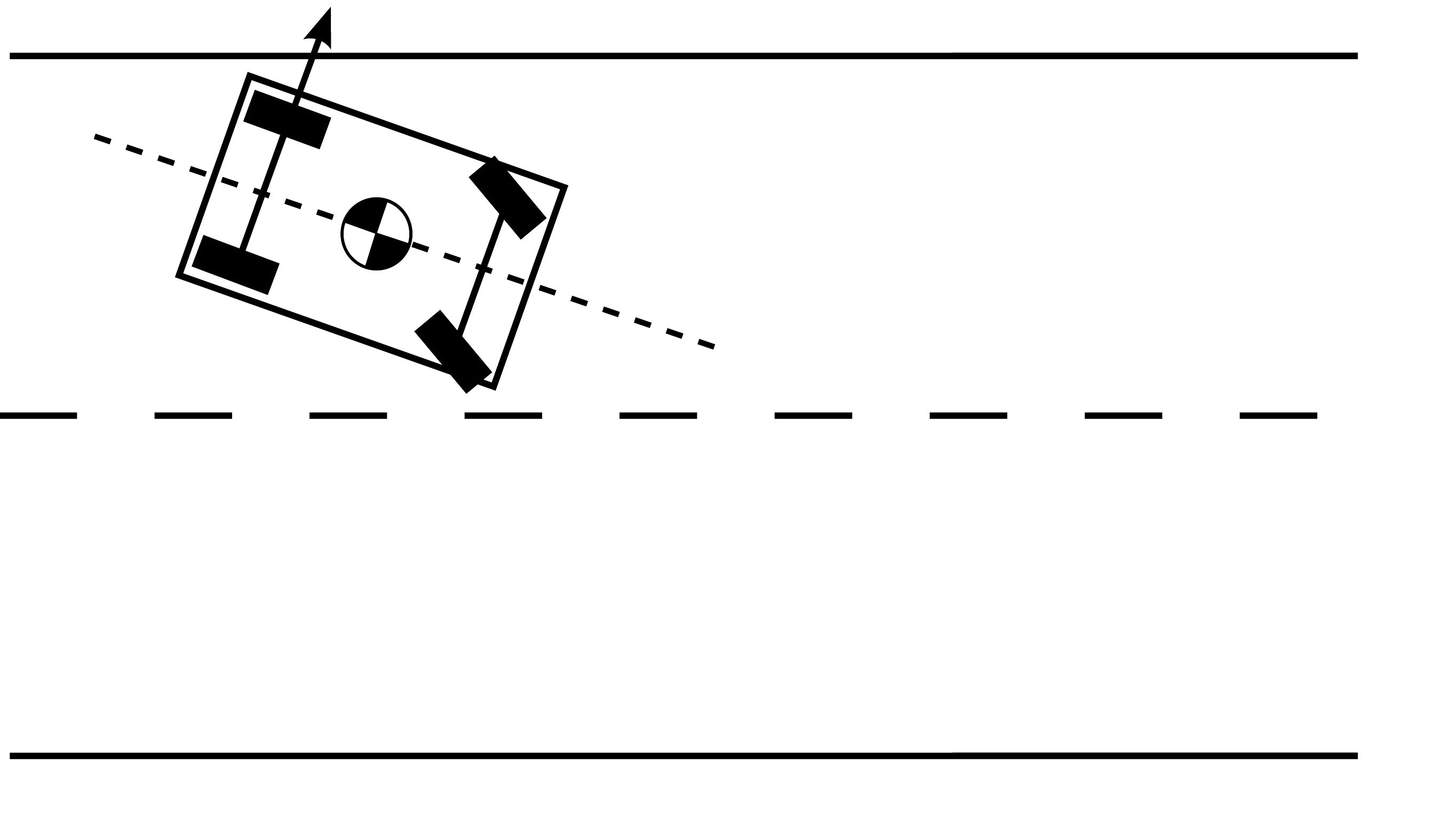
        \caption{A scenario in which the dynamical object at time $N$, i.e., $o_N$, does not cover the ego vehicle reference position at time $N$.}
        \label{fig:problemstatement}
    \end{subfigure}
    \begin{subfigure}{\columnwidth}
    \centering
        \def\svgwidth{0.75\columnwidth}
        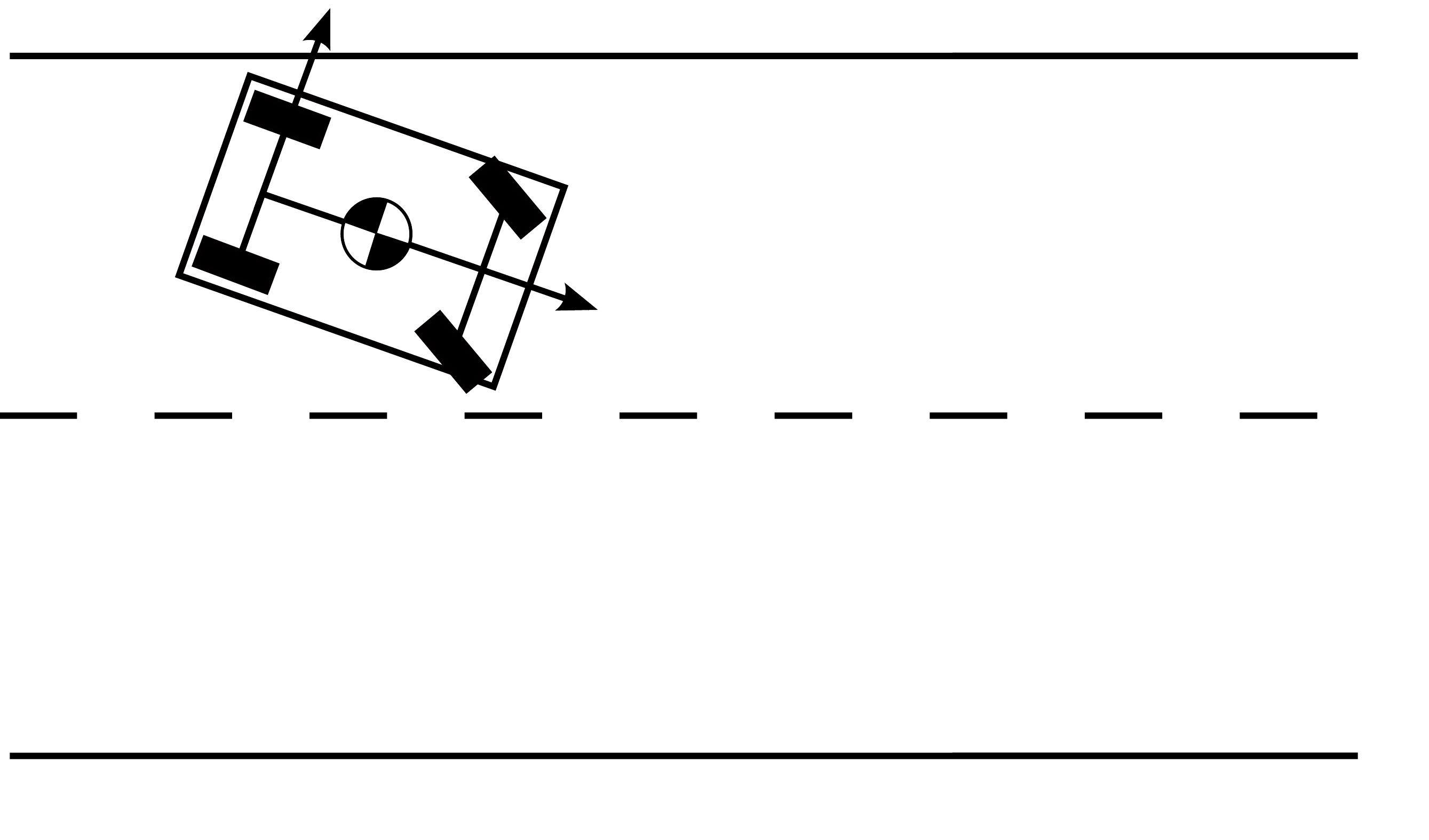
        \caption{A scenario in which the dynamical object at time $N$, i.e., $o_N$, covers the reference position at time instance $N$.}
        \label{fig:problemstatement2}
    \end{subfigure}
    \caption{Two illustrations introduce the problem statement and show the typical trade-off between driving towards a goal versus minimal-risk driving.}
\end{figure}

In various on-road situations (e.g., the lane change motion tasks shown in Fig.~\ref{fig:problemstatement}) the automated vehicle needs to plan a safe, comfortable and drivable motion. To do so, the automated vehicle is expected to obey the structured environment, (i.e., the road infrastructure, traffic rules, etc.), interact with mixed-traffic (manually driven vehicles, cyclists, pedestrians), and account for their future uncertain movements in various situations (e.g. bad light, poor road markings, bad weather, dense traffic). 

The vehicle, with coordinate system $(x,y)$, is expected to drive from an initial state given by $(x_0,y_0,\theta_0,v_0)$, located on top of the coordinate system (i.e., $x_0,y_0,\theta_0=0$) , towards the reference state $r_N$ which is positioned relative to the initial condition of the vehicle and contains a planar position and longitudinal velocity $v$ as a reference, i.e., $(x_N,y_N,\theta_N,v_{N})$.  Here, $N$ denotes the horizon of the trajectory generator to be planned towards the goal state, i.e., the look-ahead time. Throughout this work, we assume that the reference state, $\mathbf{r}_N$, is located in the center of the rightmost lane on the road at a distance of $N\cdot t_s\cdot \bar{v}$, where $t_s$ represents the sampling time between each trajectory element and $\bar{v}$ represents the maximally allowed velocity, provided by the infrastructure. Furthermore, it is assumed that this reference state is known a-priori, and could be provided by a global planner (e.g., an $\text{A}^*$ algorithm~\cite{Likhachev2009933}).  The road infrastructure is assumed in this work to be described by an integer set of polynomials, representing the boundaries of the lanes, i.e.:
\begin{align}
    \label{eq:poly}
    \tilde{y}_{i} = c_0^{(i)}+c_1^{(i)}\tilde{x}_{i}+c_2^{(i)}\tilde{x}_{i}^2+c_3^{(i)}\tilde{x}_{i}^3
\end{align}
where $i\in\mathbb{Z}$, the superscript $^{(i)}$ selects the polynomial constants of the $i$-th line measurement, i.e., $c_j^{(i)},\:j\in[0,1,2,3]$, and ($\tilde{x},\tilde{y}$) denote the planar coordinates of the $i$-th line relative to the vehicle initial condition. Note, that polynomial functions of arbitrary order can be used following our methodology. Finally, the vehicle is expected to safely drive along with the infrastructure in the presence of objects (e.g., vulnerable road users (VRUs) or other vehicles). The state of the objects $\mathbf{o}_k\forall k\in[0\hdots N]$ is denoted as $(x^o_k,y^o_k,\theta^o_k)$, which is measured and predicted with respect to the initial condition of the vehicle. To include the predicted motion of the object in motion planning, the state of the object is assumed to be exactly measurable/predictable up to time $N$.

The trajectory of the vehicle is expected to minimize risks by abiding by the traffic rules while respecting the presence of other road users. Furthermore, the planner is expected to plan for a period that allows anticipating the environment (i.e., tens of seconds down to 3-4 seconds~\cite{Calvert202042}). Note that minimizing risks while driving towards a reference $\mathbf{r}_N$ could result in an optimal trajectory that does not reach the reference. To illustrate this phenomenon, we reconsider Fig.~\ref{fig:problemstatement} in a different scenario depicted in Fig.~\ref{fig:problemstatement2}. In this scenario, the vehicle is expected to avoid risk by, e.g., steering the vehicle back into the center of the lane while avoiding a dangerous situation with the object, meaning it will inherently not reach its reference state $\mathbf{r}_N$.

Next to the planar movement of the vehicle (i.e., $x,y,\theta$), a notion of time needs to be included to be able to incorporate the risk of moving objects. As a result, the planar movement is augmented with the longitudinal velocity of the vehicle, denoted by $v$. The main challenge of a trajectory generator is, hence, to find a reference trajectory, $\mathcal{T}$, for the vehicle to follow, which takes into account all the above aspects,
\begin{equation}
   \mathcal{T}= (x_k,y_k,\theta_k,v_{k})\:\forall k\in[0\hdots N]
\end{equation}
where the trajectory length is $N$ and satisfies kinematic vehicle behavior and comfort-driven constraints while minimizing the risk to the passengers and surrounding road users. 

\section{Model-predictive planning for risk-averse driving}~\label{sec:method}
In this work, we approach the trajectory planning problem from a model-predictive control (MPC) perspective. The characteristics of an MPC are in line with our requirements, namely, minimizing risk (i.e., formulating the cost as a function of risk), while satisfying kinematic vehicle dynamical behavior and constraints (e.g., limited steering input or longitudinal acceleration).  
\subsection{Formulation of the vehicle model}
\begin{figure}[h!]
    \centering
    \def\svgwidth{0.60\columnwidth}
\begingroup%
  \makeatletter%
  \providecommand\color[2][]{%
    \errmessage{(Inkscape) Color is used for the text in Inkscape, but the package 'color.sty' is not loaded}%
    \renewcommand\color[2][]{}%
  }%
  \providecommand\transparent[1]{%
    \errmessage{(Inkscape) Transparency is used (non-zero) for the text in Inkscape, but the package 'transparent.sty' is not loaded}%
    \renewcommand\transparent[1]{}%
  }%
  \providecommand\rotatebox[2]{#2}%
  \newcommand*\fsize{\dimexpr\f@size pt\relax}%
  \newcommand*\lineheight[1]{\fontsize{\fsize}{#1\fsize}\selectfont}%
  \ifx\svgwidth\undefined%
    \setlength{\unitlength}{379.84251969bp}%
    \ifx\svgscale\undefined%
      \relax%
    \else%
      \setlength{\unitlength}{\unitlength * \real{\svgscale}}%
    \fi%
  \else%
    \setlength{\unitlength}{\svgwidth}%
  \fi%
  \global\let\svgwidth\undefined%
  \global\let\svgscale\undefined%
  \makeatother%
  \begin{picture}(1,0.89552239)%
    \lineheight{1}%
    \setlength\tabcolsep{0pt}%
    \put(0,0){\includegraphics[width=\unitlength,page=1]{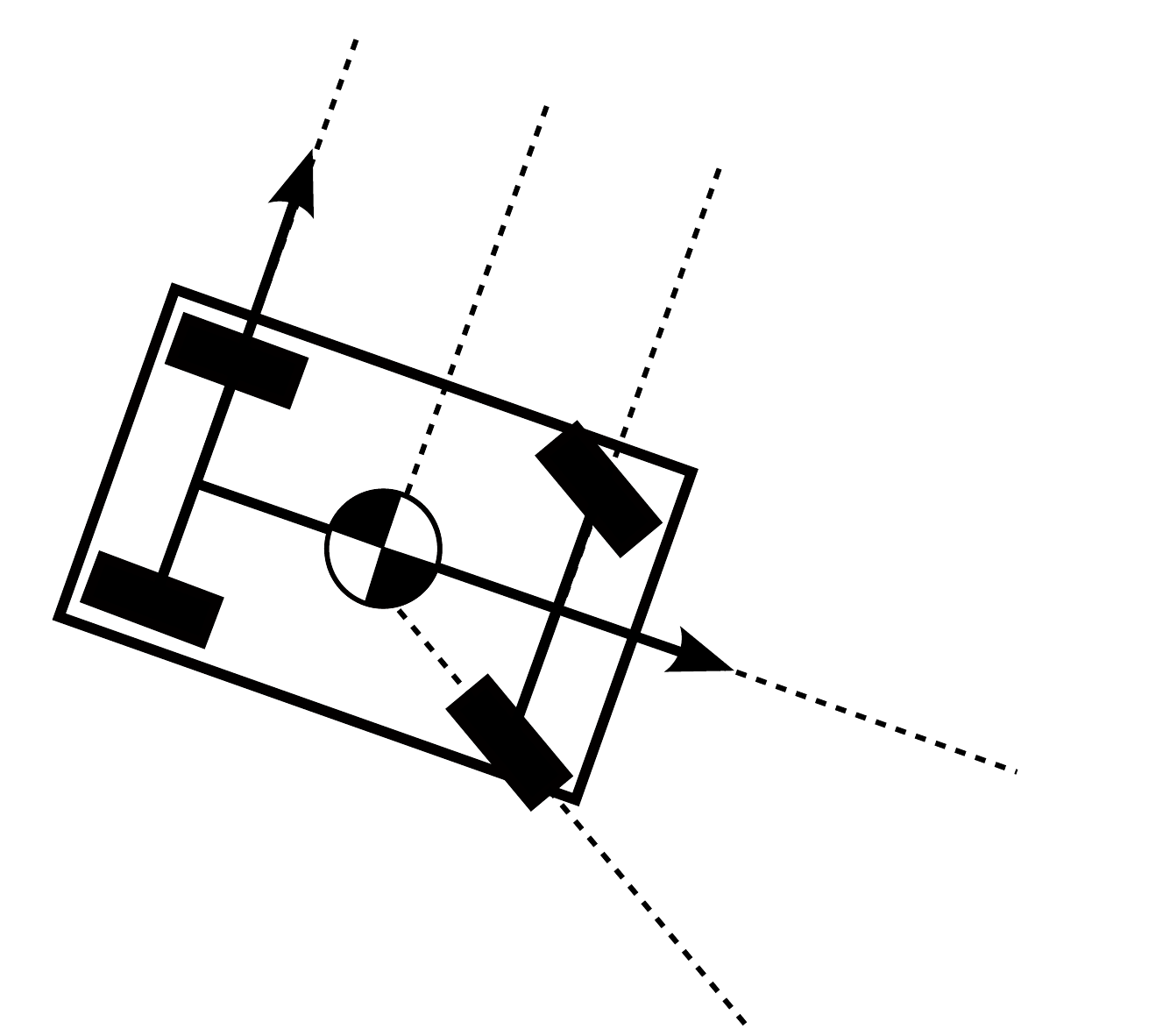}}%
    \put(0.61006759,0.35268583){\color[rgb]{0,0,0}\rotatebox{-1.098766}{\makebox(0,0)[lt]{\lineheight{1.25}\smash{\begin{tabular}[t]{l}$x,v$\end{tabular}}}}}%
    \put(0.17104855,0.67326789){\color[rgb]{0,0,0}\rotatebox{-0.552413}{\makebox(0,0)[lt]{\lineheight{1.25}\smash{\begin{tabular}[t]{l}$y$\end{tabular}}}}}%
    \put(0,0){\includegraphics[width=\unitlength,page=2]{kinematic_model_standalone.pdf}}%
    \put(0.448214,0.82740204){\color[rgb]{0,0,0}\makebox(0,0)[lt]{\lineheight{1.25}\smash{\begin{tabular}[t]{l}$L$\end{tabular}}}}%
    \put(0.3402279,0.78126749){\color[rgb]{0,0,0}\makebox(0,0)[lt]{\lineheight{1.25}\smash{\begin{tabular}[t]{l}$L_r$\end{tabular}}}}%
    \put(0.48387505,0.6924149){\color[rgb]{0,0,0}\makebox(0,0)[lt]{\lineheight{1.25}\smash{\begin{tabular}[t]{l}$L_f$\end{tabular}}}}%
    \put(0,0){\includegraphics[width=\unitlength,page=3]{kinematic_model_standalone.pdf}}%
    \put(0.72958431,0.05240956){\color[rgb]{0,0,0}\makebox(0,0)[lt]{\lineheight{1.25}\smash{\begin{tabular}[t]{l}$\delta$\end{tabular}}}}%
  \end{picture}%
\endgroup%

    \caption{Kinematic bicycle model, with the planar reference frame in the center of the rear axis.}
    \label{fig:3_model}
\end{figure}

The model, used for MPC synthesis, considered in this work is the non-linear kinematic bicycle model. This model is particularly suitable for larger sampling times, i.e., larger time horizons while maintaining numerical stability~\cite{9366415}. The model is depicted in Fig.~\ref{fig:3_model} and is described in discrete time, using forward Euler discretization with sampling time $t_s$. We denote the state of the model as $\mathbf{x}=[x,y,\theta,v]$ and its inputs as $\mathbf{u}=[\delta,a]$, and the state-update rules are denoted as follows
\begin{align}
\begin{aligned}
    x_{k+1} =& x_k+t_sv_k\cos{\theta_k},\\
    y_{k+1} =& y_k+t_sv_k\sin{\theta_k},\\
    \theta_{k+1} =& \theta_k+t_s\frac{v_k}{L}\tan{\delta_k},\\
    v_{k+1} =& v_k+t_sa_k.
\end{aligned}\label{eq:vehiclemodel}
\end{align}
where the states $x,y,\theta$ denote the planar movement relative to the initial condition, as denoted in Fig.~\ref{fig:problemstatement}, $v$ represents the velocity in the direction of the $x$-axis, $\delta$ denotes the steering angle of the front axle, $a$ denotes the acceleration in the direction of the $x$-axis and finally $L$ denotes the wheelbase of the vehicle. The state vector $\mathbf{x}$ and the input vector $\mathbf{u}$ belong to the sets
\begin{subequations}
\label{eq:boundedsets}
\begin{align}
    \mathbf{X}=&\{\mathbf{x}\in\mathbb{R}^4|\underline{\mathbf{x}}\leq\mathbf{x}\leq\overline{\mathbf{x}}\},\\ \mathbf{U}=&\{\mathbf{u}\in\mathbb{R}^2|\underline{\mathbf{u}}\leq\mathbf{u}\leq\overline{\mathbf{u}}\},
\end{align}
\end{subequations}
respectively. The bounds of the sets are governed by infrastructural, dynamic, and comfort limitations.
\subsection{Formulation of the risk}
\begin{figure}[h!]
    \centering
    \begin{subfigure}{\columnwidth}
        \includegraphics[width=0.9\columnwidth]{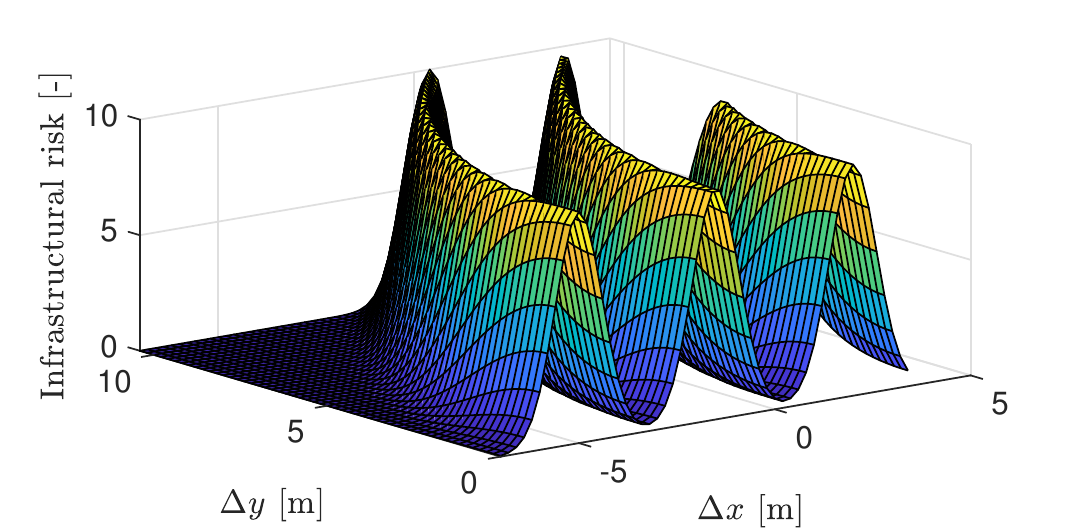}
        \caption{Example infrastructural risk-field of a slightly curved two-lane highway.}
        \label{fig:infra_pfield}
    \end{subfigure}
    \begin{subfigure}{\columnwidth}
        \includegraphics[width=0.9\columnwidth]{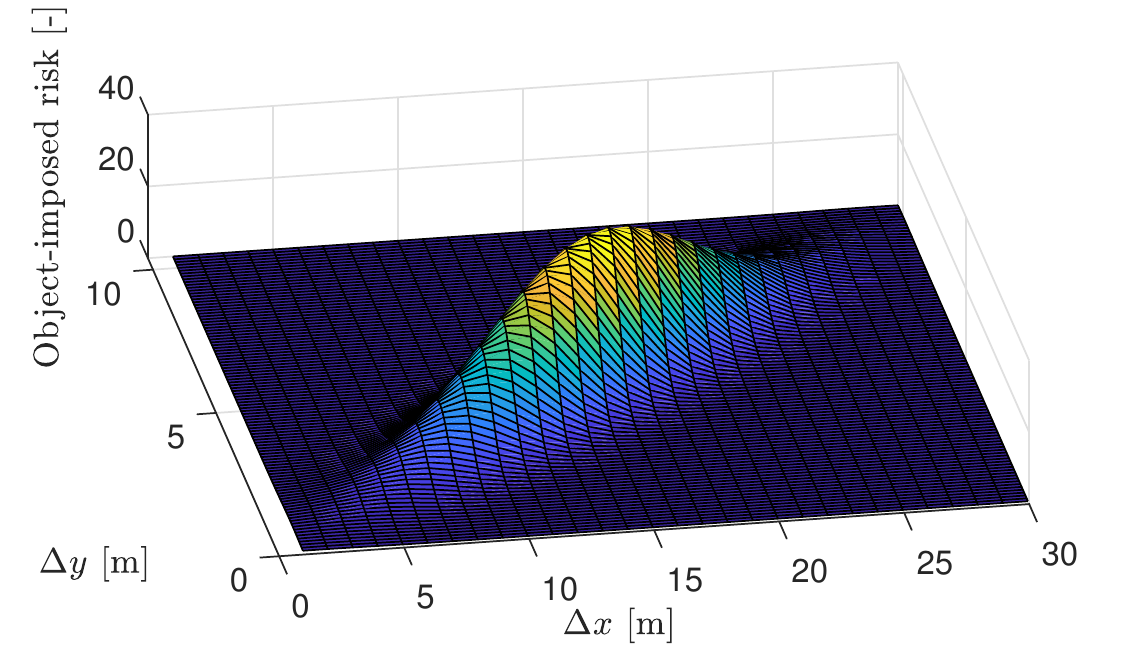}
        \caption{Example object risk-field.}
        \label{fig:object_pfield}
    \end{subfigure}
    \caption{Two illustrative examples of the proposed risk-field definitions. $\Delta x$ denotes the $x$-distance from the vehicle's initial position and $\Delta y$ denotes the $y$-distance from the vehicle's initial position.}
\end{figure}

Similar to~\cite{Li2022122}, we consider two sources of risk, namely, infrastructural risk and object-imposed risk. The infrastructural risk is described as the safety risk taken to travel across certain road boundaries ~\cite{Li2022122}. We define the infrastructural risk, at time instance $k$, by $\mathcal{U}^{I}_k$. Furthermore, we define a risk-free situation, i.e., $\mathcal{U}^{I}_k\approx0$ if and only if the vehicle is driving in the center of a lane and assuming that all lanes are driveable. The risk-function of the infrastructure, perceived by the vehicle at time instance $k$, can be described as a set of structured Gaussian functions, using the polynomial road description~\eqref{eq:poly}, as follows
\begin{align}
    \label{eq:infra}
    \mathcal{U}^{I}_k=&\sum_{i=1}^{N_d}A_Ie^{-\frac{\left(y_k-\tilde{y}_k^{(i)}(\tilde{x}^{(i)}_k)\right)^2}{2 \sigma^2}},\forall i\in \mathbb{Z}^+
\end{align}
where $N_d$ represents the number of road lines, $A_I$ represents the amplitude of the risk and $\sigma$ represents the standard deviation of the risk. The amplitude and standard deviation of the infrastructural risk are assumed constant, though may depend on the type of lanes that are divided (e.g., a lane with oncoming traffic may be riskier to drive on, hence the amplitude may be higher). This risk-field description is a slight adaption of the description in~\cite{9206298}, where only a two-lane highway is considered. An illustrational example of an infrastructural potential field is provided in Fig.~\ref{fig:infra_pfield}.

The relative risk of dynamical objects is described using a repulsive potential field formulated as a multivariate, potentially rotated, Gaussian distribution as follows:
\begin{align}
    \mathcal{U}^{O}_k=&\sum_{n=1}^{N_o}A_Oe^{-\frac{f_o^{(n)}(\mathbf{x}_k,\mathbf{o}_k)}{2}}\label{eq:objects}\\
    f_o^{(n)}=&\begin{bmatrix}x_k-\widehat{x}^{(n)}_k& y_k-\widehat{y}^{(n)}_k\end{bmatrix}R\Sigma^{-1} R^\intercal\begin{bmatrix}x_k-\widehat{x}^{(n)}_k\\ y_k-\widehat{y}^{(n)}_k\end{bmatrix}\nonumber\\
    \Sigma=&\begin{bmatrix}{\sigma_x^2}&0\\0&{\sigma_y^2}\end{bmatrix},\quad R=\begin{bmatrix}\cos\widehat{\theta}^{(n)}_k&-\sin\widehat{\theta}^{(n)}_k\\\sin\widehat{\theta}^{(n)}_k&\cos\widehat{\theta}^{(n)}_k\end{bmatrix}\nonumber
\end{align}
where $(\widehat{x}_k^{(n)},\widehat{y}_k^{(n)},\widehat{\theta}_k^{(n)})$ represents the state of the $n$-th object, $N_o$ represents the number of objects, $A_O\in\mathbb{R}^N_d$ represents for each object the maximal risk, $\sigma_x$ represents the longitudinal standard deviation of the perceived risk of the object and $\sigma_y$ represents the lateral standard deviation of the perceived risk of the object. Note, again, that for simplicity the amplitude and standard deviations are assumed constant, though they can be varied based on the type of object (e.g., pedestrian versus a car) or for example the velocity at which the ego vehicle is driving. An illustrational example of such a risk field is provided in Fig.~\ref{fig:object_pfield}.
\subsection{Formulation of the risk-averse trajectory generator}
Using the definition of the vehicle dynamic model~\eqref{eq:vehiclemodel}, its limitations~\eqref{eq:boundedsets} and the definition of the risk field for the infrastructure~\eqref{eq:infra} and objects~\eqref{eq:objects}, the optimization problem for a risk-averse trajectory generator can be composed.
The cost and the constraints of a risk-averse trajectory generator are formulated by the following model-predictive program:
\begin{subequations}
\label{eq:mpcprogram}
    \begin{align}
        J^\star_{0\rightarrow N}(\mathbf{x}_s)=&\min_{u,x}\sum_{k=1}^{N-1}\ell_\mathbf{x}\left(\mathbf{x}_k-\mathbf{r}_k\right)+\sum_{k=0}^{N-1}\ell_\mathbf{u}\left(\mathbf{u}_k\right)\nonumber\\+m&(\mathbf{x}_{N}-\mathbf{r}_{N})+\mathcal{U}_k^O+\mathcal{U}_k^I\label{eq:pathmpc1}\\\text{s.t.}\quad\mathbf{x}_0=&\mathbf{x}_s\label{eq:pathmpc2}\\
        \quad \mathbf{x}_{k+1}=&f(\mathbf{x}_k,\mathbf{u}_k),\quad\forall k\geq 0\label{eq:pathmpc3}\\
        \quad g(\mathbf{x}_k,\mathbf{u}_k)\leq&\: 0,\quad \forall k\geq 0\label{eq:pathmpc4}
    \end{align}
\end{subequations}
where, in~\eqref{eq:pathmpc1}, the functions $\ell_{\mathbf{x}}(\cdot),\:\ell_{\mathbf{u}}(\cdot)$ quadratically penalizes the control error $\mathbf{x}_k-\mathbf{r}_k$, and the input $\mathbf{u}_k$, respectively. Furthermore, the term $m(\cdot)$ quadratically penalizes the terminal state at time $N$. The state $\mathbf{x}_s$ in~\eqref{eq:pathmpc2} represents the state at present time. The function $f$ in~\eqref{eq:pathmpc3} represents the non-linear dynamics as depicted in~\eqref{eq:vehiclemodel}. Finally, the constraints in~\eqref{eq:pathmpc4} enforce the boundedness of the state $\mathbf{x}$ and the input $\mathbf{u}$ through enforcement of the sets~\eqref{eq:boundedsets}. Note, that the stage cost $\ell_\mathbf{x}$ is not used in this work (i.e., equal to zero), due to the fact that there is only the requirement to drive to a goal at time instance $N$, hence, only the terminal cost $m(\cdot)$ is used to penalize the state. 
\section{Simulation study}~\label{sec:simulation}
The non-linear program is implemented in a real-time operating optimization framework CasADi~\cite{Andersson2018} using the IPOPT Interior-point embedded solver~\cite{Wachter200625}. The trajectory generator is deployed in four different studies, where urban- and highway driving in a risk-defined environment forms the basis. Throughout the first three simulations, the same vehicle parameter set is used to show the versatility of the algorithm without having to tune it for every separate use case. For the urban case, slight adjustments are made. The parameter-set is provided in Table~\ref{table:pars}, where it has to be noted, that any parameters which are not denoted are either $0$ in the case of stage weights or not set in the case of state bounds. In each of the test cases, the computational time is denoted\footnote{The results are generated on an Intel Core i7-10850 2.7GHz platform.}. The optimization problem is initiated at each instance $t_s$, hence, by real-time we mean that the MPC trajectory generator can generate a new solution within a period of $t_s$ seconds. Finally, in the simulation, we assume perfect tracking of the trajectory by the vehicle model.
\begin{table}
\small
\begin{tabular}{p{1.2cm} p{3.8cm} p{1.4cm} p{0.5cm}}\toprule
\textbf{Parameter}              & \textbf{Description} & \textbf{Value} & \textbf{Unit}\\
\toprule
$t_s$ & Sampling time & 0.75 & s\\
$N$            & MPC horizon & 10 & -\\
$\ell_\mathbf{u}$           & Input weight & $[1,100]$ & -\\
$m$&Terminal weight&$[1,0.01,0,0]$&-\\
$A_I$&Infrastructure risk amplitude&100&-\\
$A_O$&Object risk amplitude&1000&-\\
$\sigma$&Infrastructure risk deviation&1.3&$\text{m}$\\
$\sigma_x$&Longitudinal risk deviation&20&$\text{m}$\\
$\sigma_y$&Lateral risk deviation&1.3&$\text{m}$\\
$\underline{a},\overline{a}$&Vehicle acceleration bounds&$[-4,0.5]$&$\text{ms}^{-2}$\\
$\underline{\delta},\overline{\delta}$&Vehicle steering angle bounds&$[-0.1,0.1]$&$\text{rad}$\\
$\underline{y},\overline{y}$&Lateral displacement bounds&$[1,9.5]$&$\text{m}$\\
$\underline{v},\overline{v}$&Vehicle velocity bounds&$[0,10]$&$\text{ms}^{-1}$\\
\toprule
\end{tabular}
\caption{Simulation parameters used for case \RNum{1} up to \RNum{3} and partially for case \RNum{4}.}
\label{table:pars}
\end{table}
\subsection{Case \RNum{1}: single overtake}
\begin{figure}[h!]
    \centering
    \begin{subfigure}{\columnwidth}
        \includegraphics[width=1\columnwidth]{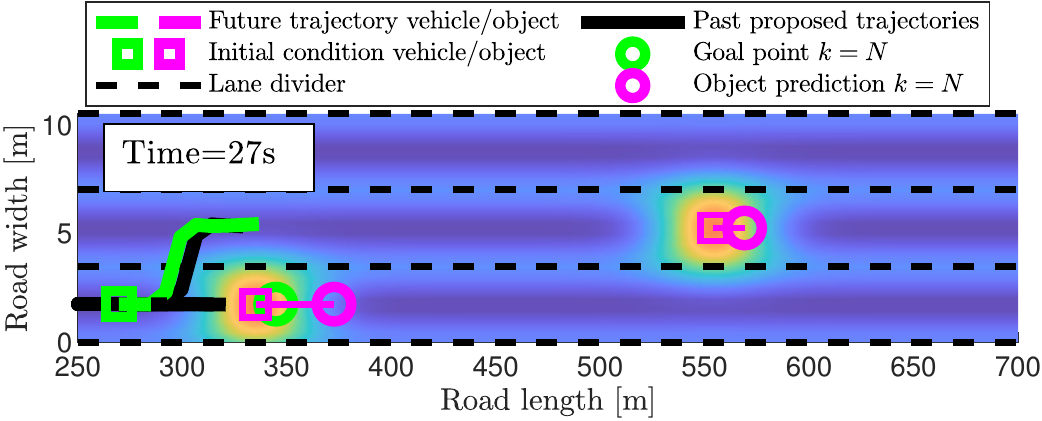}
    \end{subfigure}
    \begin{subfigure}{\columnwidth}
        \includegraphics[width=1\columnwidth]{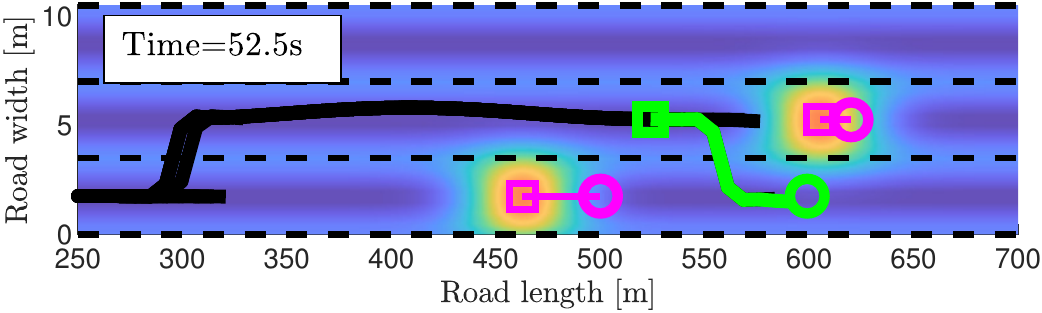}
    \end{subfigure}
    \begin{subfigure}{\columnwidth}
        \includegraphics[width=1\columnwidth]{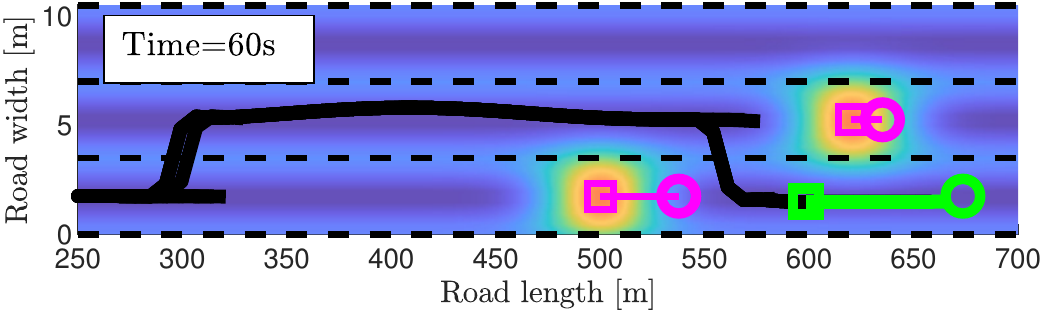}
    \end{subfigure}
    \caption{Simulation results of case \RNum{1} with two vehicles driving slowly on the two rightmost lanes, motivating a single lane change.}
    \label{fig:case1}
\end{figure}

In the first case, we consider a scenario where the vehicle drives in the rightmost lane and approaches several slow-driving vehicles. The objects are driving at a velocity of $5\text{ms}^{-1}$ (bottom lane) and $2\text{ms}^{-1}$ (middle lane), respectively as typical in highway traffic jam situations. The results are shown in Fig.~\ref{fig:case1} at several interesting points in time. First, the vehicle approaches the slow-driving vehicle in the bottom lane and overtakes the vehicle. Then, the vehicle has three choices, it can either pass the object on the top or bottom lane, or it can slow down and drive behind the vehicle. In this case, it is of minimal cost (i.e., a combination of risk and benefit) to return to the original lane, such that the vehicle can, once again, drive towards the goal point. The average computational time of the planner in this scenario is $30.4\text{ms}$ with a standard deviation of $5.1\text{ms}$.
\subsection{Case \RNum{2}: double overtake}
\begin{figure}[h!]
    \centering
    \begin{subfigure}{\columnwidth}
        \includegraphics[width=1\columnwidth]{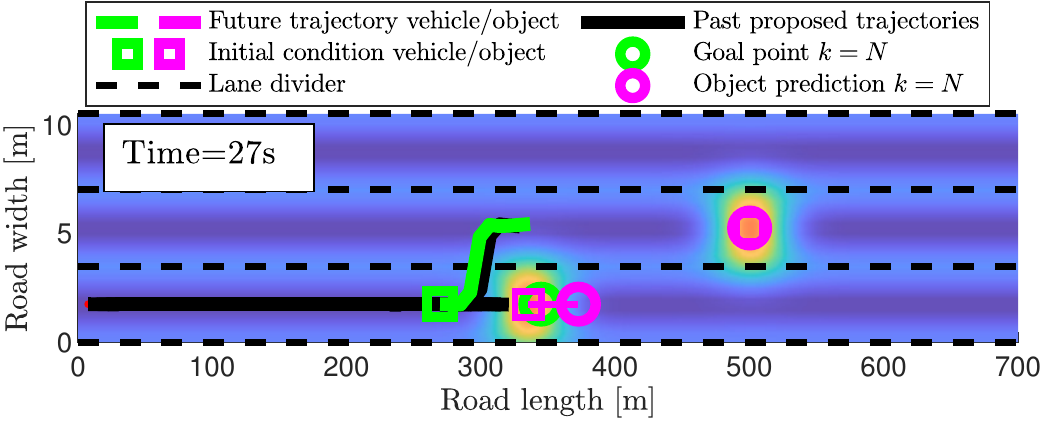}
    \end{subfigure}
    \begin{subfigure}{\columnwidth}
        \includegraphics[width=1\columnwidth]{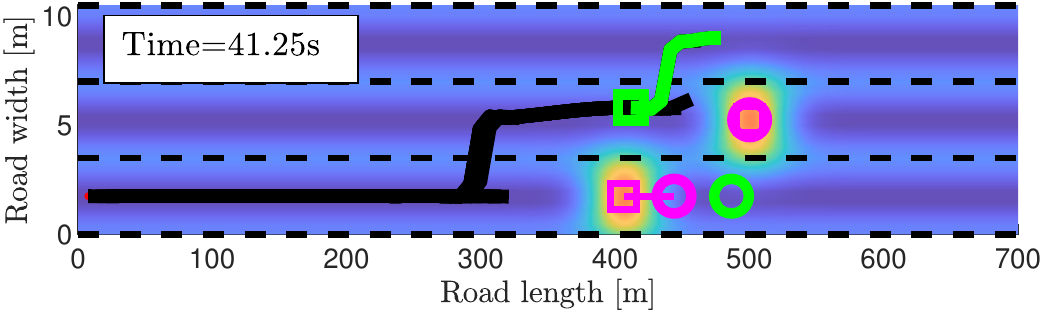}
    \end{subfigure}
    \begin{subfigure}{\columnwidth}
        \includegraphics[width=1\columnwidth]{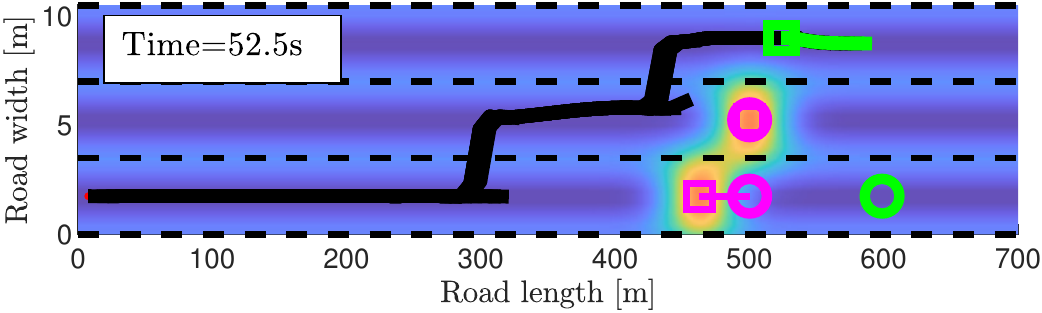}
    \end{subfigure}
    \caption{Simulation results of case \RNum{2} with two vehicles driving slowly on the two rightmost lanes, motivating a double lane change.}
    \label{fig:case2}
\end{figure}

In the second case, we consider the scenario where the vehicle starts on the rightmost lane, overtakes a slow vehicle, and then decides whether to slow down, or to once again overtake due to a static vehicle in the middle lane. The objects are driving at a velocity of $5\text{ms}^{-1}$ (bottom lane) and $0\text{ms}^{-1}$ (middle lane), respectively. The results are shown in Fig.~\ref{fig:case2}.  The first overtake maneuver is identical to case I. Once the vehicle approaches the static vehicle, it can be observed that a lane change to the right would compromise the taken risk due to the proximity of the right-lane object prediction at horizon $N$ ahead. As a result, the MPC takes the anticipatory decision to take another lane change to the left, at the cost of a higher lateral distance error to the desired goal point. The average computational time of the planner in this scenario is $30.2\text{ms}$ with a standard deviation of $4.1\text{ms}$.
\subsection{Case \RNum{3}: blocked lanes, slowly driving participants}
\begin{figure}[h!]
    \centering
    \begin{subfigure}{\columnwidth}
        \includegraphics[width=1\columnwidth]{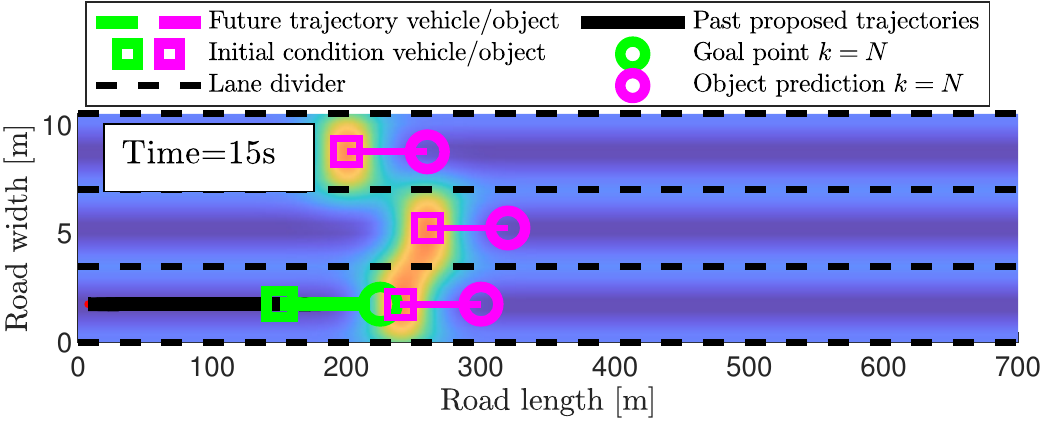}
    \end{subfigure}
    \begin{subfigure}{\columnwidth}
        \includegraphics[width=1\columnwidth]{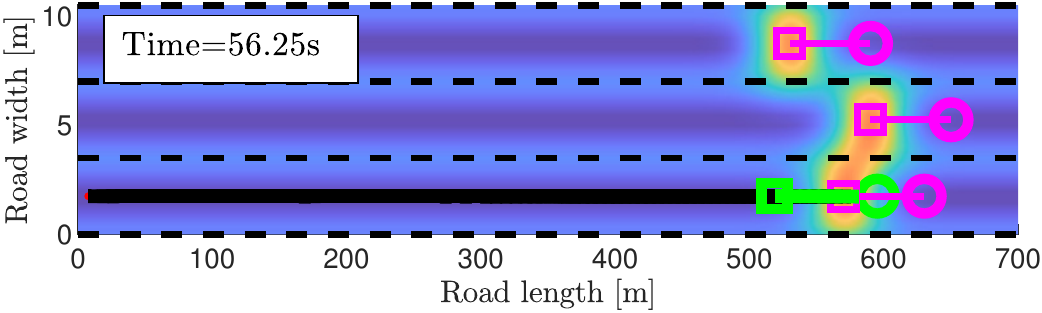}
    \end{subfigure}
    \begin{subfigure}{\columnwidth}
    \centering
        \includegraphics[width=\columnwidth]{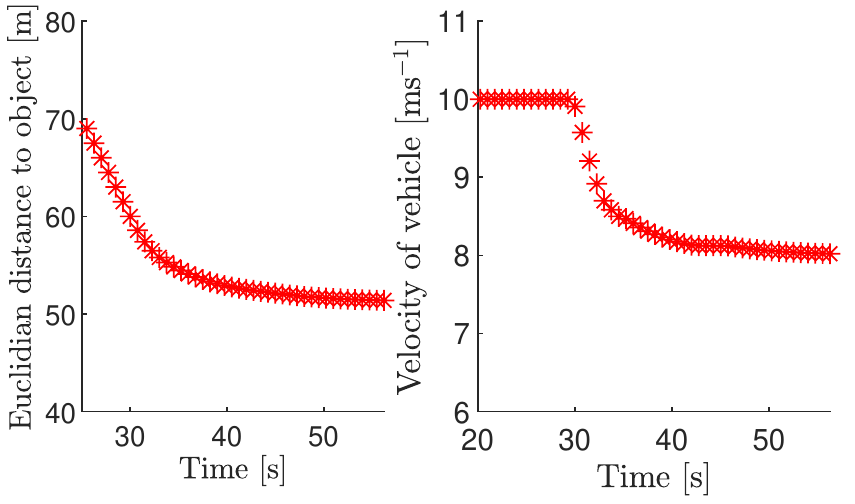}
    \end{subfigure}
    \caption{Simulation results of case \RNum{3} with three vehicles slowly driving and blocking the lanes.}
    \label{fig:case3}
\end{figure}

In the third case, we consider the scenario where the three-lane road is blocked by slowly driving objects. The objects are all driving at a velocity of $8\text{ms}^{-1}$ and their longitudinal position initial condition deviates to challenge the algorithm. The results are shown in Fig.~\ref{fig:case3}. It can be observed that the MPC notices that there is only one real sensible solution, which is to stay inside the lane and brake to the same velocity as the succeeding vehicles. Note, that the steady-state euclidean distance between the vehicles can be tuned using the longitudinal parameter $\sigma_x$.
The average computational time of the planner in this scenario is $28.2\text{ms}$ with a standard deviation of $1.6\text{ms}$.
\subsection{Case \RNum{4}: VRU crossing with potential late object detection}
\begin{figure}[h!]
    \includegraphics[width=1\columnwidth]{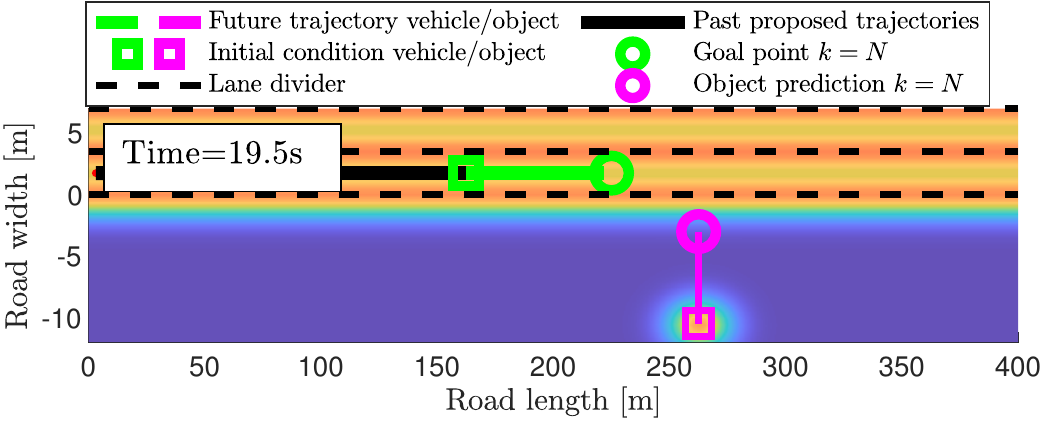}
    \caption{Initial setup for the VRU simulation case \RNum{4}.}
    \label{fig:UC4_0}
\end{figure}

\begin{figure}[h!]
        \centering
        \begin{subfigure}{\columnwidth}
            \includegraphics[width=1\columnwidth]{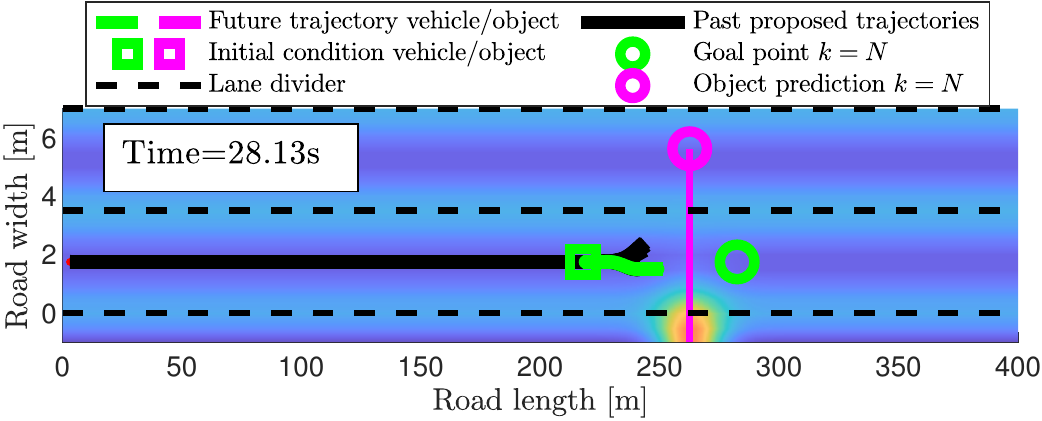}
        \end{subfigure}
        \begin{subfigure}{\columnwidth}
            \includegraphics[width=1\columnwidth]{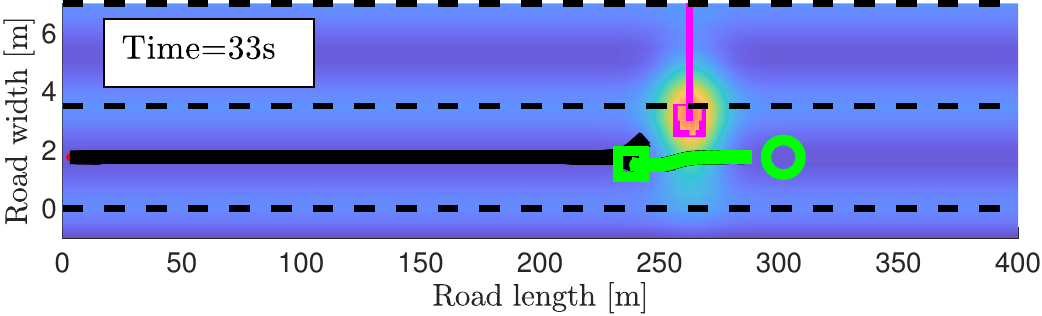}
        \end{subfigure}
        \begin{subfigure}{\columnwidth}
            \includegraphics[width=\columnwidth]{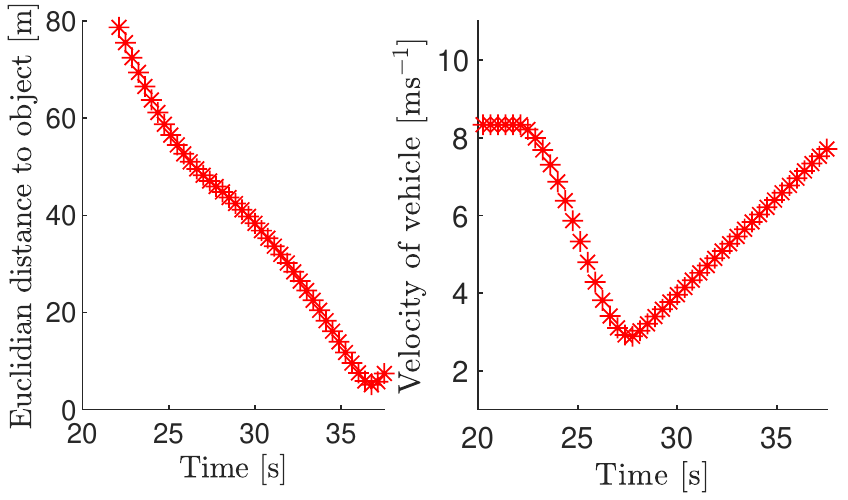}
        \end{subfigure}
        \caption{Simulation results of case \RNum{4} with timely object detection.}
        \label{fig:case4detection}
    \end{figure}

\begin{figure}
    \centering
    \begin{subfigure}{\columnwidth}
        \includegraphics[width=1\columnwidth]{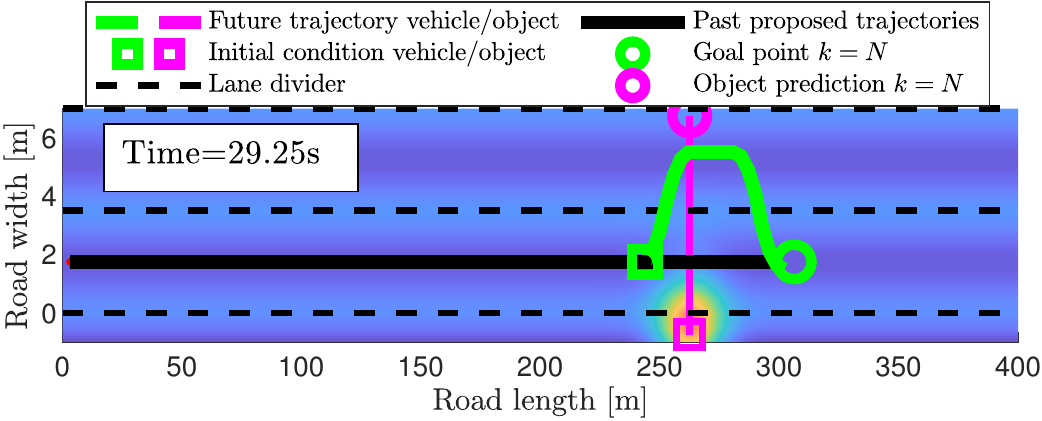}
    \end{subfigure}
    \begin{subfigure}{\columnwidth}
        \includegraphics[width=1\columnwidth]{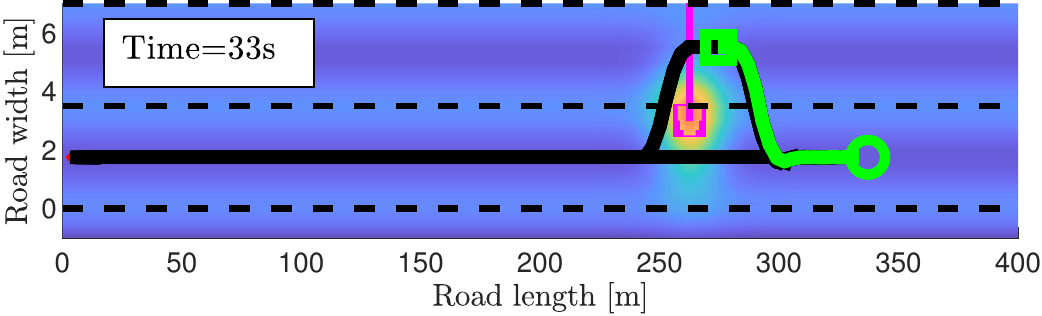}
    \end{subfigure}
    \begin{subfigure}{\columnwidth}
            \includegraphics[width=\columnwidth]{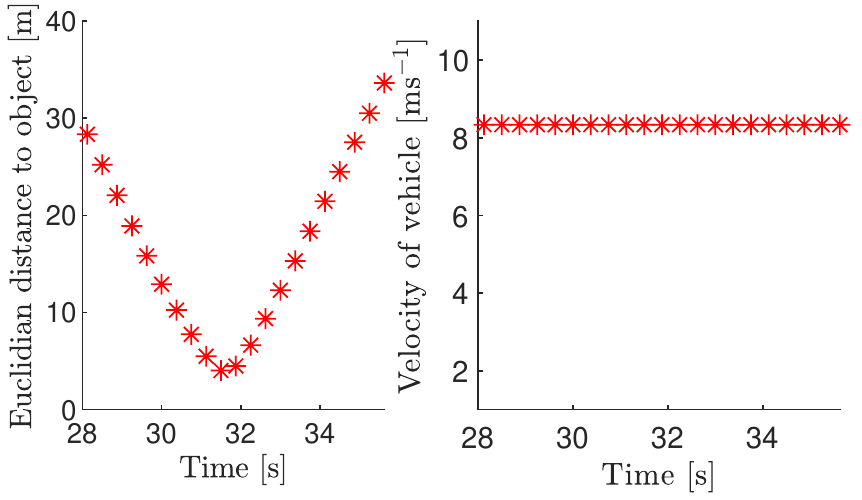}
        \end{subfigure}
    \caption{Simulation results of case \RNum{4} with late object detection.}
    \label{fig:case4nodetection}
\end{figure}

 Whereas case \RNum{1} to \RNum{3} presented highway maneuvers, the fourth and final case will cover an urban driving setting, showing the applicability of the MPC trajectory generator for design domains other than highway driving. In this case, we consider the vehicle driving along a two-lane road. On this two-lane road, a VRU is about to cross the road as depicted in Fig.~\ref{fig:UC4_0}. This case is deemed to be highly relevant and largely responsible for VRU incidents (up to $51\%$~\cite{Yue2020119}). The proposed simulations are twofold. First, we study how the MPC trajectory generator acts when it sees the VRU and can anticipate its actions based on these observations. Subsequently, we study how the trajectory generator acts when the VRU observation is obstructed and the detection of the VRU is late (advocated in~\cite{Yue2020119} as a relevant root cause in this particular scenario). Transitioning from a highway scenario towards an urban scenario requires an adaptation of MPC parameters. The revised parameters, different from Table~\ref{table:pars}, are proposed in Table~\ref{table:pars2}. In an urban scenario, the maneuvers likely require a higher tracking fidelity, for this reason, the sampling time is reduced and, as a result, the horizon is increased to achieve a similar look-ahead time concerning the previous cases. Moreover, the reduced longitudinal velocity results in a change of bound and a change in the longitudinal standard deviation of the object risk-field. Finally, an urban scenario with two-lane roads usually allows bi-directional traffic. As a result, the infrastructural potential field amplitude is increased~(as shown in Table~\ref{table:pars} and also in Fig.~\ref{fig:UC4_0} by the infrastructural risk field color). The VRU crosses the road at a velocity of $1\text{ms}^{-1}$ and the ego vehicle starts at an initial velocity of $8.33\text{ms}^{-1}$. First, we consider the results without late VRU detection, as depicted in Fig.~\ref{fig:case4detection}. It is clear that, from $23\text{s}$ onward the vehicle starts to anticipate and brake, down to a velocity of $3\text{ms}^{-1}$, to wait for the VRU to pass. Once the VRU has passed, the vehicle accelerates back to its original velocity, which is a natural and intuitive expected behavior of the vehicle. Now, consider a late detection of the VRU (coming, e.g., from bad weather conditions) of which the results are depicted in Fig.~\ref{fig:case4nodetection}. The detection occurs at around $30\text{s}$ (i.e., the top graph), a time at which there is only $13\text{m}$ of distance left concerning the object. At this time, the trajectory generator path indicates that braking is no longer an option and an evasive maneuver towards the left lane is necessary. Once the evasive maneuver is finished, the vehicle returns to its lane while a safe distance has been kept with the VRU. The average computational time of the planner for the early VRU detection is $43.5\text{ms}$ with a standard deviation of $4.2\text{ms}$. The average computational time of the planner in the late detection case is $48.2\text{ms}$ with a standard deviation of $8.1\text{ms}$. The increased mean computational time can be explained largely by the fact that, due to the higher required fidelity of the path, the horizon is doubled.
\begin{table}
\small
\begin{tabular}{p{1.2cm} p{3.8cm} p{1.4cm} p{0.5cm}}\toprule
\textbf{Parameter}              & \textbf{Description} & \textbf{Value} & \textbf{Unit}\\
\toprule
$t_s$ & Sampling time & 0.38 & s\\
$N$            & MPC horizon & 20 & -\\
$A_I$&Infrastructure risk amplitude&200&-\\
$\sigma_x$&Longitudinal risk deviation&10&$\text{m}$\\
$\underline{y},\overline{y}$&Lateral displacement bounds&$[1,6]$&$\text{m}$\\
$\underline{v},\overline{v}$&Vehicle velocity bounds&$[0,8.33]$&$\text{ms}^{-1}$\\
\toprule
\end{tabular}
\caption{Changed parameters w.r.t. Table~\ref{table:pars} for case \RNum{4}.}
\label{table:pars2}
\end{table}
 

\section{Conclusion}~\label{sec:conclusion}
In this work, we have proposed a monolithic model-predictive trajectory generator that, as opposed to many works in literature, does not use a set of pre-defined decisions or trajectories to decide on actions to progress towards a goal point while minimizing risk. This offers a simpler solution to trajectory generation, which can be used for longer horizons while being computationally feasible. The proposed method explicitly considers predicted trajectories of other road users and it's shown effective in various use cases. The method is shown to be effective by demonstrating it in four relevant scenarios in which the approach shows the ability to make logical and safe decisions in both highway and urban scenarios. Furthermore, the planner is challenged in a more difficult urban scenario, where the late detection of a pedestrian (for example happening due to bad weather) is compared with the timely detection. The results show that this approach can function safely under both circumstances. Future work includes the adaptation of the MPC problem formulation and the potential fields to incorporate even more complex situations and be robust for measurement uncertainty, induced by predictions or faults (for example, the lack of fault isolability in object fault detection~\cite{9655077}). Moreover, adaptive potential fields will be studied to make the trajectory generator suitable for miscellaneous environments.
\section{Acknowledgments}
This work is supported by the EU Horizon 2020 R{\&}D program under grant agreement No. 861570, project SAFE-UP (proactive SAFEty systems and tools for a constantly UPgrading road environment).
    \printbibliography

@ARTICLE{9206298,
    author={Siebenrock, F. and Gunther, M. and Hohmann, S.},
    title={LTV-MPC based trajectory planning considering uncertain object prediction through adaptive potential fields},
    journal={IEEE Conference on Control Technology and Applications},
    year={2020},
    pages={666-672},
}

@ARTICLE{9366415,
  author={Laurense, Vincent A. and Gerdes, J. Christian},
  journal={IEEE Transactions on Control Systems Technology}, 
  title={Long-Horizon Vehicle Motion Planning and Control Through Serially Cascaded Model Complexity}, 
  year={2022},
  volume={30},
  number={1},
  pages={166-179}}

@Article{Andersson2018,
  Author = {Joel A E Andersson and Joris Gillis and Greg Horn
            and James B Rawlings and Moritz Diehl},
  Title = {{CasADi} -- {A} software framework for nonlinear optimization
           and optimal control},
  Journal = {Mathematical Programming Computation},
  Year = {2018},
}

@ARTICLE{Wachter200625,
author={Wächter, A. and Biegler, L.T.},
title={On the implementation of an interior-point filter line-search algorithm for large-scale nonlinear programming},
journal={Mathematical Programming},
year={2006},
volume={106},
number={1},
pages={25-57}
}

@ARTICLE{DeFreitas2021,
author={De Freitas, J. and Censi, A. and Smith, B.W. and Lillo, L.D. and Anthony, S.E. and Frazzoli, E.},
title={From driverless dilemmas to more practical commonsense tests for automated vehicles},
journal={Proceedings of the National Academy of Sciences of the United States of America},
year={2021},
volume={118},
number={11}
}

@ARTICLE{Khatib1985500,
author={Khatib, O.},
title={Real-time obstacle avoidance for manipulators and mobile robots},
journal={IEEE International Conference on Robotics and Automation},
year={1985},
pages={500-505}
}

@ARTICLE{González20161135,
author={González, D. and Pérez, J. and Milanés, V. and Nashashibi, F.},
title={A Review of Motion Planning Techniques for Automated Vehicles},
journal={IEEE Transactions on Intelligent Transportation Systems},
year={2016},
volume={17},
number={4},
pages={1135-1145}
}

@ARTICLE{Likhachev2009933,
author={Likhachev, M. and Ferguson, D.},
title={Planning long dynamically feasible maneuvers for autonomous vehicles},
journal={International Journal of Robotics Research},
year={2009},
volume={28},
number={8},
pages={933-945}
}

@ARTICLE{Ji2017952,
author={Ji, J. and Khajepour, A. and Melek, W.W., Sr. and Huang, Y.},
title={Path planning and tracking for vehicle collision avoidance based on model predictive control with multiconstraints},
journal={IEEE Transactions on Vehicular Technology},
year={2017},
volume={66},
number={2},
pages={952-964}
}

@ARTICLE{Gutjahr20171586,
author={Gutjahr, B. and Gröll, L. and Werling, M.},
title={Lateral Vehicle Trajectory Optimization Using Constrained Linear Time-Varying MPC},
journal={IEEE Transactions on Intelligent Transportation Systems},
year={2017},
volume={18},
number={6},
pages={1586-1595}
}

@ARTICLE{Dixit20181061,
author={Dixit, S. and Montanaro, U. and Fallah, S. and Dianati, M. and Oxtoby, D. and Mizutani, T. and Mouzakitis, A.},
title={Trajectory Planning for Autonomous High-Speed Overtaking using MPC with Terminal Set Constraints},
journal={IEEE Conference on Intelligent Transportation Systems},
year={2018},
pages={1061-1068},
}

@ARTICLE{Calvert202042,
author={Calvert, S.C. and Schakel, W.J. and van Lint, J.W.C.},
title={A generic multi-scale framework for microscopic traffic simulation part II – Anticipation Reliance as compensation mechanism for potential task overload},
journal={Transportation Research Part B: Methodological},
year={2020},
volume={140},
pages={42-63}
}

@ARTICLE{Papadoulis201912,
author={Papadoulis, A. and Quddus, M. and Imprialou, M.},
title={Evaluating the safety impact of connected and autonomous vehicles on motorways},
journal={Accident Analysis and Prevention},
year={2019},
volume={124},
pages={12-22}
}

@ARTICLE{Huang2016232,
author={Huang, Z. and Wu, Q. and Ma, J. and Fan, S.},
title={An APF and MPC combined collaborative driving controller using vehicular communication technologies},
journal={Chaos, Solitons and Fractals},
year={2016},
volume={89},
pages={232-242},
}

@ARTICLE{Hang202014458,
author={Hang, P. and Lv, C. and Huang, C. and Cai, J. and Hu, Z. and Xing, Y.},
title={An Integrated Framework of Decision Making and Motion Planning for Autonomous Vehicles Considering Social Behaviors},
journal={IEEE Transactions on Vehicular Technology},
year={2020},
volume={69},
number={12},
pages={14458-14469}
}

@ARTICLE{Yue2020119,
author={Yue, L. and Abdel-Aty, M. and Wu, Y. and Zheng, O. and Yuan, J.},
title={In-depth approach for identifying crash causation patterns and its implications for pedestrian crash prevention},
journal={Journal of Safety Research},
year={2020},
volume={73},
pages={119-132},
}

@ARTICLE{Li2022122,
author={Li, L. and Gan, J. and Ji, X. and Qu, X. and Ran, B.},
title={Dynamic Driving Risk Potential Field Model under the Connected and Automated Vehicles Environment and Its Application in Car-Following Modeling},
journal={IEEE Transactions Intelligent Transportation Systems},
year={2022},
volume={23},
number={1},
pages={122-141},
}

@INPROCEEDINGS{9655077,
  author={van der Ploeg, Chris and Smit, Robin and Siagkris-Lekkos, Alexis and Benders, Frank and Silvas, Emilia},
  booktitle={European Control Conference}, 
  title={Anomaly Detection from Cyber Threats via Infrastructure to Automated Vehicle}, 
  year={2021},
  volume={},
  number={},
  pages={1788-1794}}
\end{document}